\documentclass[11pt]{amsart}
\usepackage{amsmath,amssymb}
\newtheorem{theorem}{Theorem}
\newtheorem{corollary}[theorem]{Corollary}

\newtheorem{proposition}[theorem]{Proposition}

\newtheorem*{definition}{Definition}

\begin{document}

\title[Conformal Scalar Curvature Equation]{On the Conformal Scalar Curvature Equation and Related Problems}
\author{Simon Brendle}

\maketitle 

\section{The Yamabe problem}

A fundamental result in two-dimensional Riemannian geometry is the uniformization theorem, which assers that any 
Riemannian metric on a compact surface is conformally related to a metric of constant curvature. There are two ways of generalizing this result to higher dimensions. One problem, proposed by E.~Calabi, is concerned with the existence of 
K\"ahler-Einstein metrics (or, more generally, extremal K\"ahler metrics). The other one was proposed in 1960 by H.~Yamabe \cite{Yamabe}. Yamabe claimed the following theorem: 

\begin{theorem}
\label{yamabe.problem}
Let $M$ be a compact manifold of dimension $n \geq 3$ without boundary, and let 
$g$ be a Riemannian metric on $M$. Then there exists a metric 
$\tilde{g}$ which is conformally related to $g$ and has constant scalar curvature.
\end{theorem}

However, Yamabe's proof contained a serious gap. Theorem \ref{yamabe.problem} was first
verified by N.~Trudinger \cite{Trudinger} under the additional assumption that $g$ has non-positive scalar
curvature. T.~Aubin \cite{Aubin1} proved Theorem \ref{yamabe.problem} assuming that
$n \geq 6$ and $(M,g)$ is not locally conformally flat. The remaining cases were
solved by R.~Schoen \cite{Schoen1} using the positive mass theorem (see \cite{Schoen-Yau1}, \cite{Witten}).
A.~Bahri \cite{Bahri} gave an alternative proof in the locally conformally flat case.

It is well known that the Yamabe problem is equivalent to the solvability of a semilinear elliptic PDE. Indeed, if $\tilde{g}$ is conformally related to $g$, then we can write $\tilde{g} = u^{\frac{4}{n-2}} \, g$, where $u$ is a positive function on $M$. The scalar curvature of $\tilde{g}$ is related to the scalar curvature of $g$ by the formula 
\begin{equation} 
\label{transformation.law.for.scalar.curvature}
R_{\tilde{g}} \, u^{\frac{n+2}{n-2}} = R_g \, u - \frac{4(n-1)}{n-2} \, \Delta_g u.
\end{equation} 
Here, $\Delta_g$ denotes the Laplace operator associated with the metric $g$. In particular, the metric $\tilde{g} = u^{\frac{4}{n-2}} \, g$ has scalar curvature $c$ if and only if $u$ is a positive solution of the PDE 
\begin{equation} 
\label{yamabe.pde} 
\frac{4(n-1)}{n-2} \, \Delta_g u - R_g \, u + c \, u^{\frac{n+2}{n-2}} = 0. 
\end{equation} 
The Yamabe problem admits a natural variational characterization. A function $u$ satisfies (\ref{yamabe.pde}) for some 
$c \in \mathbb{R}$ if and only if $u$ is a critical point of the Yamabe functional 
\begin{equation} 
\label{yamabe.energy}
E_g(u) = \frac{\int_M \big ( \frac{4(n-1)}{n-2} \, |du|_g^2 + R_g \, u^2 \big ) \, dvol_g}{\big ( \int_M u^{\frac{2n}{n-2}} \, dvol_g \big )^{\frac{n-2}{n}}}. 
\end{equation}
The functional $E_g$ resembles the Sobolev quotient associated with the embedding of $W^{1,2}(M,g)$ into $L^{\frac{2n}{n-2}}(M,g)$. Moreover, the functional $E_g$ is closely related to the total scalar curvature functional. Indeed, it follows from (\ref{transformation.law.for.scalar.curvature}) that $E_g(u) = \mathcal{E}(u^{\frac{4}{n-2}} \, g)$, where $\mathcal{E}$ denotes the normalized Einstein-Hilbert action: 
\[\mathcal{E}(g) = \frac{\int_M R_g \, dvol_g}{\text{\rm Vol}_g(M)^{\frac{n-2}{n}}}.\] 
Finally, the Yamabe constant of $g$ is defined as the infimum of the Yamabe functional $E_g(u)$ over all positive functions on $M$: 
\begin{equation} 
\label{yamabe.constant}
Y(M,g) = \inf_{0 < u \in C^1(M)} E_g(u). 
\end{equation}
The key step in the proof of Theorem \ref{yamabe.problem} is an upper bound, due to T.~Aubin and R.~Schoen, for the Yamabe constant $Y(M,g)$. Once this inequality is known, the existence of a minimizer of $E_g$ follows from standard PDE arguments.

\begin{theorem} 
Assume that $(M,g)$ is not conformally equivalent to the round sphere $S^n$. Then $Y(M,g) < Y(S^n)$, where $Y(S^n)$ denotes the Yamabe energy of the round metric on $S^n$.
\end{theorem}

The Yamabe PDE (\ref{yamabe.pde}) has a unique solution if $Y(M,g) \leq 0$. However, uniqueness fails for $Y(M,g) > 0$. As an example, consider the sphere $S^n$ equipped with the round metric. This example is very special, because the Yamabe PDE is invariant under the conformal transformations on $S^n$. By a result of M.~Obata \cite{Obata}, every metric of constant scalar curvature which is conformal to the round metric on $S^n$ is given by the pull-back of the round metric under a conformal diffeomorphism. Hence, if $g$ is the round metric on $S^n$, then every solution to (\ref{yamabe.pde}) is minimizing, and the space of solutions can be identified with the unit ball $B^{n+1}$.

As another example, let $g$ be the product metric on $S^{n-1}(1) \times S^1(L)$. In this case, all solutions of the Yamabe PDE are rotationally symmetric. If the length $L$ of the $S^1$-factor is sufficiently small, then the Yamabe PDE has a unique solution (which is constant). On the other hand, the Yamabe PDE has many non-minimizing solutions if $L$ is large. We refer to \cite{Schoen2} for a detailed discussion of this example.

In this paper, we will discuss two PDE problems that arise in connection with the Yamabe 
problem. The first problem is concerned with the set of all constant scalar curvature metrics 
in a given conformal class. In particular, we are interested in knowing whether or not this set 
is compact. A compactness result for solutions of (\ref{yamabe.pde}) would be useful in 
developing a Morse theory for the Yamabe functional (see \cite{Schoen3} and \cite{Schoen4} for 
details). 

The second problem deals with the gradient flow for the Yamabe functional. 
This leads to the following evolution equation for the Riemannian metric: 
\begin{equation} 
\label{yamabe.flow}
\frac{\partial}{\partial t} g(t) = -(R_{g(t)} - r_{g(t)}) \, g(t), \qquad g(0) = g_0. 
\end{equation}
Here, $R_{g(t)}$ denotes the scalar curvature of $g(t)$ and $r_{g(t)}$ is 
the mean value of the scalar curvature of $g(t)$. We are interested in the longtime behavior of this evolution equation: 
since (\ref{yamabe.flow}) is the gradient flow to the Yamabe functional, one might expect that the 
flow converges to a metric of constant scalar curvature as $t \to \infty$.

\section{Compactness of the set of constant scalar curvature metrics in a given conformal class}

Fix a compact Riemannian manifold $(M,g)$ of dimension $n \geq 3$ with positive Yamabe constant. In this section, we will address the question of whether the set of solutions to the Yamabe PDE 
\begin{equation} 
\tag{$\star$} 
\frac{4(n-1)}{n-2} \, \Delta_g u - R_g \, u + 4n(n-1) \, u^{\frac{n+2}{n-2}} = 0 
\end{equation} 
is compact in the $C^2$-topology. In some applications, it is convenient to approximate ($\star$) by subcritical problems of the form 
\begin{equation} 
\tag{$\star_\delta$} 
\frac{4(n-1)}{n-2} \, \Delta_g u - R_g \, u + 4n(n-1) \, u^{\frac{n+2}{n-2} - \delta} = 0, 
\end{equation} 
where $\delta \in [0,\frac{4}{n-2})$. Given any real number $\delta_0 \in [0,\frac{4}{n-2})$, we denote by $\Phi(M,g;\delta_0)$ the set of all positive functions that satisfy either the Yamabe PDE ($\star$) or one of the subcritical problems: 
\begin{align*} 
&\Phi(M,g;\delta_0) \\ 
&= \Big \{ 0 < u \in C^2(M): \text{$u$ is a solution of ($\star_\delta$) for some $\delta \in [0,\delta_0]$} \Big \}. 
\end{align*} 
Clearly, compactness fails if $(M,g)$ is conformally equivalent to the round sphere. The 
following result provides sufficient conditions for the set $\Phi(M,g;\delta_0)$ to 
be compact.

\begin{theorem} 
\label{compactness}
Assume that $(M,g)$ is not conformally equivalent to the round sphere $S^n$. Moreover, assume that one of the following conditions holds:
\begin{itemize} 
\item[(i)] $(M,g)$ is locally conformally flat 
\item[(ii)] $3 \leq n \leq 7$ 
\item[(iii)] $n \geq 8$ and $|W_g(p)| + |\nabla W_g(p)| > 0$ for all $p \in M$ 
\end{itemize} 
Given any real number $\delta_0 \in [0,\frac{4}{n-2})$, there exists a constant $\Lambda$, depending only on $\delta_0$ and the background metric $g$, such that $\|u\|_{C^2(M)} \leq \Lambda$ for all functions $u \in \Phi(M,g;\delta_0)$.
\end{theorem}

R.~Schoen proved Theorem \ref{compactness} in the locally conformally flat case (see \cite{Schoen2}, \cite{Schoen3}). Moreover, Schoen outlined a strategy for proving the theorem in the non-locally conformally flat case. In \cite{Li-Zhu}, Y.Y.~Li and M.~Zhu followed this strategy to prove the theorem in dimension $3$. O.~Druet \cite{Druet} proved compactness in dimensions $4$ and $5$.

The case $n \geq 6$ is more subtle, and requires a careful analysis of the local properties of the background metric $g$ near a blow-up point. A well-known conjecture asserts that the Weyl tensor should vanish to an order greater than $[\frac{n-6}{2}]$ at a blow-up point (see \cite{Schoen4}). This conjecture is known as the Weyl Vanishing Conjecture. It has been verified in dimensions $6$ and $7$ by F.~Marques and, independently, by Y.Y.~Li and L.~Zhang \cite{Li-Zhang1}. Using this result and the positive mass theorem, these authors were able to prove compactness for all $n \leq 7$. In dimension $n \geq 8$, Li and Zhang proved that any blow-up point $p$ satisfies $W_g(p) = 0$ and $\nabla W_g(p) = 0$ (see \cite{Li-Zhang1}). In particular, this proves the Weyl Vanishing Conjecture in dimensions $8$ and $9$. Moreover, this result shows that blow-up cannot occur if $|W_g(p)| + |\nabla W_g(p)| > 0$ for all $p \in M$.

The preceeding discussion raises the question as to whether the Weyl Vanishing Conjecture holds in dimension $n \geq 10$. Y.Y.~Li and L.~Zhang \cite{Li-Zhang2} proved the Weyl Vanishing Conjecture in dimensions $10$ and $11$. Very recently, M.~Khuri, F.~Marques, and R.~Schoen \cite{Khuri-Marques-Schoen} were able to verify the Weyl Vanishing Conjecture up to dimension $24$. This result, combined with the positive mass theorem, can be used to rule out blow-up for $n \leq 24$ (see \cite{Khuri-Marques-Schoen}, Theorem 1.1): 

\begin{theorem} 
\label{compactness.24} 
Assume that $(M,g)$ is not conformally equivalent to the round sphere $S^n$. Moreover, assume that $n \leq 24$ and $M \setminus \{p\}$ is a spin manifold. For every real number $\delta_0 \in [0,\frac{4}{n-2})$, there exists a constant $\Lambda$, depending only on $\delta_0$ and the background metric $g$, such that $\|u\|_{C^2(M)} \leq \Lambda$ for all functions $u \in \Phi(M,g;\delta_0)$.
\end{theorem}

A key ingredient in the proof of both Theorem \ref{compactness} and Theorem \ref{compactness.24} is the Pohozaev identity. Let $\Omega$ be an open subset of $M$ with smooth boundary, and let $\nu$ be the outward-pointing unit normal vector field along $\partial \Omega$. Moreover, suppose that $V$ is a vector field defined on $\Omega$. 
Using the divergence theorem, we obtain 
\begin{align*} 
&- \int_\Omega \frac{1}{2} \, \Big [ \nabla^i V^j +\nabla^j V^i - \frac{2}{n} \, \text{\rm div}_g \, V \, g^{ij} \Big ] \, \partial_i u \, \partial_j u \, dvol_g \\ 
&- \frac{n-2}{4(n-1)} \int_\Omega R_g \, u \, \Big [ \langle V,\nabla u \rangle + \frac{n-2}{2n} \, u \, \text{\rm div}_g \, V \Big ] \, dvol_g \\ 
&+ \int_{\partial \Omega} \langle V,\nabla u \rangle \, \langle \nabla u,\nu \rangle - \int_{\partial \Omega} \frac{1}{2} \, |\nabla u|^2 \, \langle V,\nu \rangle \\ 
&+ \frac{n-2}{2n} \int_{\partial \Omega} u \, \langle \nabla u,\nu \rangle \, \text{\rm div}_g \, V + \frac{(n-2)^2}{2} \, \int_{\partial \Omega} u^{\frac{2n}{n-2}-\delta} \, \langle V,\nu \rangle \\ 
&= \frac{n-2}{2n} \int_\Omega u \, \langle \nabla u,\nabla(\text{\rm div}_g \, V) \rangle \, dvol_g \\ 
&- \frac{(n-2)^2}{2} \, \delta \, \int_\Omega u^{\frac{n+2}{n-2}-\delta} \, \langle V,\nabla u \rangle \, dvol_g  
\end{align*}
whenever $u$ is a solution of ($\star_\delta$). This relation can be viewed as a generalization 
of the classical Pohozaev identity in Euclidean space. The proof of Theorem \ref{compactness} 
rests on a precise analysis of the various terms in the Pohozaev identity. For example, in low 
dimensions and in the locally conformally flat case, one can use the positive mass theorem to 
show that blow-up is inconsistent with the Pohozaev identity.
We refer to \cite{Druet-Hebey2} and \cite{Schoen4} for excellent surveys on this subject. 

Theorem \ref{compactness} can be used to obtain information concerning the number of constant 
scalar curvature metrics. This idea was proposed by R.~Schoen in \cite{Schoen4}. Assume that the 
assumptions of Theorem \ref{compactness} are satisfied. Moreover, assume that all positive 
solutions to ($\star$) are non-degenerate. This implies that the number of solutions to ($\star$) is 
finite. Let $N_k$ be the number of solutions to ($\star$) with Morse index $k$. 
By compactness, every solution to ($\star_\delta$) is close (in the $C^2$-topology) to a solution of ($\star$) if $\delta > 0$ is sufficiently small. 
Conversely, near each solution of ($\star$) there is exactly one solution to ($\star_\delta$), and this 
solution has the same Morse index as the original one. Consequently, the subcritical problem ($\star_\delta$) has exactly $N_k$ solutions with Morse index $k$. 
Since the functional associated with ($\star_\delta$) satisfies the Palais-Smale condition, the Morse 
inequalities hold for the subcritical problem ($\star_\delta$). This allows us to draw the following conclusion: 

\begin{corollary} 
\label{morse.inequalities}
Assume that the assumptions of Theorem \ref{compactness} are satisfied, and 
all solutions to ($\star$) are non-degenerate. Then 
\[\sum_{k=0}^l (-1)^{l-k} \, N_k \geq (-1)^l\] 
for $l = 0,1,2,\hdots$.
\end{corollary}

\section{Non-compactness and non-uniqueness results}

In the previous section, we have discussed various conditions that imply compactness of the 
set $\Phi(M,g;\delta_0)$. In this section, we will address the opposite question: is it 
possible to construct Riemannian manifolds $(M,g)$ such that the set of constant scalar 
curvature metrics in the conformal class of $g$ fails to be compact? Until recently, the only 
known examples where compactness fails involved non-smooth background metrics. The first result 
in this direction was established by A.~Ambrosetti and A.~Malchiodi \cite{Ambrosetti-Malchiodi}. 
This result was subsequently improved by M.~Berti and A.~Malchiodi \cite{Berti-Malchiodi}. The 
following theorem was established in \cite{Berti-Malchiodi}: 

\begin{theorem} 
\label{berti-malchiodi}
Let $k \geq 2$ and $n \geq 4k+3$. Then there exists a Riemannian metric $g$ of class $C^k$ on $S^n$ with the following properties: 
\begin{itemize}
\item[(i)] $g$ is not conformally flat 
\item[(ii)] There exists a sequence of positive functions $v_\nu$ such that $v_\nu$ is a solution of ($\star$) for all $\nu \in \mathbb{N}$ and $\sup_{S^n} v_\nu \to \infty$ as $\nu \to \infty$.
\end{itemize}
\end{theorem}

Moreover, we can choose the metric $g$ such that $g$ is close to the round metric in the 
$C^k$-norm. See \cite{Ambrosetti} for an excellent survey of this and related results. 
A related blow-up result was obtained recently by O.~Druet and E.~Hebey \cite{Druet-Hebey1}: 
they showed that blow-up can occur for problems of the form $Lu + c \, u^{\frac{n+2}{n-2}} 
= 0$, where $L$ is a small perturbation of the conformal Laplacian on $S^n$. 

In a recent paper \cite{Brendle4}, the author showed that in high dimensions blow-up can occur 
even if the background metric $g$ is $C^\infty$-smooth. More precisely, for every $n \geq 52$, there exists a smooth Riemannian metric $g$ on $S^n$ such that the set $\Phi(M,g;0)$ fails to be compact. 

\begin{theorem}
\label{blow.up}
Assume that $n \geq 52$. Then there exists a smooth Riemannian metric $g$ on $S^n$ and a sequence of positive functions $v_\nu \in C^\infty(S^n)$ ($\nu \in \mathbb{N}$) with the following properties: 
\begin{itemize}
\item[(i)] $g$ is not conformally flat 
\item[(ii)] The function $v_\nu$ is a solution of ($\star$) for all $\nu \in \mathbb{N}$
\item[(iii)] $E_g(v_\nu) < Y(S^n)$ for all $\nu \in \mathbb{N}$, and $E_g(v_\nu) \to Y(S^n)$ as $\nu \to \infty$ 
\item[(iv)] $\sup_{S^n} v_\nu \to \infty$ as $\nu \to \infty$
\end{itemize}
\end{theorem}

In particular, the sequence $v_\nu$ forms a singularity consisting of exactly one bubble.

In order to construct a metric on $S^n$ with these properties, we pick a multi-linear form $W: \mathbb{R}^n \times \mathbb{R}^n \times \mathbb{R}^n \times \mathbb{R}^n \to \mathbb{R}$ which satisfies all the algebraic properties of the Weyl tensor. 
We assume that some components $W_{ijkl}$ are non-zero, so that 
\[\sum_{i,j,k,l=1}^n (W_{ijkl} + W_{ilkj})^2 > 0.\] 
We then define a trace-free symmetric two-tensor $H_{ik}(x)$ on $\mathbb{R}^n$ by 
\[H_{ik}(x) = \sum_{p,q=1}^n W_{ipkq} \, x_p \, x_q.\] 
The proof of Theorem \ref{blow.up} rests on the following result, which is of some interest in itself
(see \cite{Brendle4}, Proposition 24). 

\begin{proposition}
\label{perturbation.argument}
Assume that $n \geq 52$. Moreover, let $g$ be a smooth metric on $\mathbb{R}^n$ of the form $g(x) = \exp(h(x))$, 
where $h(x)$ is a trace-free symmetric two-tensor on $\mathbb{R}^n$ such that $|h(x)| + |\partial h(x)| + |\partial^2 h(x)| \leq \alpha$ for all $x \in \mathbb{R}^n$, 
$h(x) = 0$ for $|x| \geq 1$, and 
$h_{ik}(x) = \mu \, (\lambda^2 - |x|^2) \, H_{ik}(x)$ for $|x| \leq \rho$. Assume that $\mu \leq 1$ and $\lambda \leq \rho \leq 1$. 
If $\alpha$ and $\rho^{2-n} \, \mu^{-2} \, \lambda^{n-10}$ are sufficiently small, then there exists a positive solution $v$ of ($\star$) satisfying 
\[\int_{\mathbb{R}^n} v^{\frac{2n}{n-2}} < \Big ( \frac{Y(S^n)}{4n(n-1)} \Big )^{\frac{n}{2}}\] 
and  
\[\sup_{|x| \leq \lambda} v(x) \geq c \, \lambda^{\frac{2-n}{2}}.\] 
Here, $c$ is a positive constant that depends only on $n$.
\end{proposition}

In order to prove Proposition \ref{perturbation.argument}, we consider the Hilbert space 
\[\mathcal{E} = \Big \{ w \in L^{\frac{2n}{n-2}}(\mathbb{R}^n) \cap W_{loc}^{1,2}(\mathbb{R}^n): dw \in L^2(\mathbb{R}^n) \Big \}\] 
with norm $\|w\|_{\mathcal{E}} = \|dw\|_{L^2(\mathbb{R}^n)}$. 
For every pair $(\xi,\varepsilon) \in \mathbb{R}^n \times (0,\infty)$, we define a closed subspace $\mathcal{E}_{(\xi,\varepsilon)} \subset \mathcal{E}$ by 
\begin{align*} 
\mathcal{E}_{(\xi,\varepsilon)} 
= \bigg \{ w \in \mathcal{E}: \int_{\mathbb{R}^n} \frac{\varepsilon^2 - |x - \xi|^2}{\varepsilon^2 + |x - \xi|^2} \, u_{(\xi,\varepsilon)}^{\frac{n+2}{n-2}} \, w &= 0 \quad \text{and} \\ 
\int_{\mathbb{R}^n} \frac{2\varepsilon \, (x_k - \xi_k)}{\varepsilon^2 + |x - \xi|^2} \, u_{(\xi,\varepsilon)}^{\frac{n+2}{n-2}} \, w &= 0 \quad \text{for $k = 1,\hdots, n$} \bigg \}. 
\end{align*} 
For every pair $(\xi,\varepsilon) \in \mathbb{R}^n \times (0,\infty)$, the function 
\[u_{(\xi,\varepsilon)}(x) = \Big ( \frac{\varepsilon}{\varepsilon^2 + |x - \xi|^2} \Big )^{\frac{n-2}{2}}\] 
is a solution of the PDE 
\[\Delta u_{(\xi,\varepsilon)} + n(n-2) \, u_{(\xi,\varepsilon)}^{\frac{n+2}{n-2}} = 0.\] 
Using the implicit function theorem, we can construct a nearby 
function $v_{(\xi,\varepsilon)}$ such that 
$v_{(\xi,\varepsilon)} - u_{(\xi,\varepsilon)} \in \mathcal{E}_{(\xi,\varepsilon)}$ and 
\[\int_{\mathbb{R}^n} \Big ( \langle dv_{(\xi,\varepsilon)},d\psi \rangle_g + \frac{n-2}{4(n-1)} \, R_g \, v_{(\xi,\varepsilon)} \, \psi - n(n-2) \, |v_{(\xi,\varepsilon)}|^{\frac{4}{n-2}} \, v_{(\xi,\varepsilon)} \, \psi \Big ) = 0\] for all test functions $\psi \in \mathcal{E}_{(\xi,\varepsilon)}$. The problem is then reduced to finding a critical point of the function 
\begin{align} 
\mathcal{F}_g(\xi,\varepsilon) &= \int_{\mathbb{R}^n} 
\Big ( |dv_{(\xi,\varepsilon)}|_g^2 + \frac{n-2}{4(n-1)} \, R_g \, v_{(\xi,\varepsilon)}^2 - (n-2)^2 \, |v_{(\xi,\varepsilon)}|^{\frac{2n}{n-2}} \Big ) \notag \\ 
&- 2(n-2) \, \Big ( \frac{Y(S^n)}{4n(n-1)} \Big )^{\frac{n}{2}}. 
\end{align} 
For abbreviation, let $\overline{H}_{ik}(x) = (1 - |x|^2) \, H_{ik}(x)$. We next define a function $F: \mathbb{R}^n \times (0,\infty) \to \mathbb{R}$ by 
\begin{align} 
F(\xi,\varepsilon) 
&= \int_{\mathbb{R}^n} \frac{1}{2} \sum_{i,k,l=1}^n \overline{H}_{il}(x) \, \overline{H}_{kl}(x) \, \partial_i u_{(\xi,\varepsilon)}(x) \, \partial_k u_{(\xi,\varepsilon)}(x) \notag \\ 
&- \int_{\mathbb{R}^n} \frac{n-2}{16(n-1)} \, \sum_{i,k,l=1}^n (\partial_l \overline{H}_{ik}(x))^2 \, u_{(\xi,\varepsilon)}(x)^2 \\ 
&+ \int_{\mathbb{R}^n} \sum_{i,k=1}^n \overline{H}_{ik}(x) \, \partial_i \partial_k u_{(\xi,\varepsilon)}(x) \, z_{(\xi,\varepsilon)}(x) \notag,
\end{align} 
where $z_{(\xi,\varepsilon)} \in \mathcal{E}_{(\xi,\varepsilon)}$ satisfies the relation 
\begin{align*} 
&\int_{\mathbb{R}^n} \Big ( \langle dz_{(\xi,\varepsilon)},d\psi \rangle - n(n+2) \, u_{(\xi,\varepsilon)}(x)^{\frac{4}{n-2}} \, z_{(\xi,\varepsilon)} \, \psi \Big ) \\ 
&= -\int_{\mathbb{R}^n} \sum_{i,k=1}^n \overline{H}_{ik} \, \partial_i \partial_k u_{(\xi,\varepsilon)} \, \psi 
\end{align*} 
for all test functions $\psi \in \mathcal{E}_{(\xi,\varepsilon)}$. The function $F(\xi,\varepsilon)$ is an even function of $\xi$. A direct computation yields 
\begin{align*} 
F(0,\varepsilon) &= -\frac{(n-2)(n+4)}{16n(n-1)(n+2)} \, |S^{n-1}| \, \sum_{i,j,k,l=1}^n (W_{ijkl} + W_{ilkj})^2 \\ 
&\hspace{10mm} \cdot \Big [ \frac{n-8}{n+4} \, \varepsilon^4 - 2 \, \varepsilon^6 + \frac{n+8}{n-10} \, \varepsilon^8 \Big ] \, \int_0^\infty (1+r^2)^{2-n} \, r^{n+3} \, dr. 
\end{align*} 
Hence, if $n \geq 52$, then the function $F(\xi,\varepsilon)$ has a critical point of the form $(0,\varepsilon_*)$, where $\varepsilon_*$ is a positive real number that depends only on $n$. 

Given any compact set $A \subset \mathbb{R}^n \times (0,\infty)$, there exists a constant $C$ such that 
\begin{align} 
\label{approximation}
&|\lambda^{-8} \, \mu^{-2} \, \mathcal{F}_g(\lambda\xi,\lambda\varepsilon) - F(\xi,\varepsilon)| \notag \\ 
&\leq C \, \lambda^{\frac{16}{n-2}} \, \mu^{\frac{4}{n-2}} + C \, \rho^{\frac{2-n}{2}} \, \mu^{-1} \, \lambda^{\frac{n-10}{2}} + C \, \rho^{2-n} \, \mu^{-2} \, \lambda^{n-10} 
\end{align}
for all $(\xi,\varepsilon) \in A$. Using (\ref{approximation}) one can show that the function $\mathcal{F}_g$ has a critical point near $(0,\lambda \varepsilon_*)$. This proves the existence of a highly concentrated solution to ($\star$). 

The following result is an immediate consequence of Proposition \ref{perturbation.argument}:

\begin{corollary} 
\label{counterexample.to.weyl.vanishing.theorem}
Assume that $n \geq 52$. Then there exists a sequence of smooth Riemannian metrics $g_\nu$ on $S^n$ and a sequence of positive functions $v_\nu \in C^\infty(S^n)$ ($\nu \in \mathbb{N}$) with the following properties: 
\begin{itemize}
\item[(i)] $\frac{4(n-1)}{n-2} \, \Delta_{g_\nu} v_\nu - R_{g_\nu} \, v_\nu + 4n(n-1) \, v_\nu^{\frac{n+2}{n-2}} = 0$ 
for all $\nu \in \mathbb{N}$ 
\item[(ii)] $E_{g_\nu}(v_\nu) < Y(S^n)$ for all $\nu \in \mathbb{N}$
\item[(iii)] The metrics $g_\nu$ converge smoothly to a metric $g_\infty$ as $\nu \to \infty$. 
\item[(iv)] There exists a sequence of points $p_\nu \in S^n$ ($\nu \in \mathbb{N}$) such that $v_\nu(p_\nu) = \sup_{S^n} v_\nu \to \infty$ as $\nu \to \infty$. 
\item[(v)] The points $p_\nu$ converge to a point $p_\infty \in S^n$ as $\nu \to \infty$. The point $p_\infty$ satisfies $W_{g_\infty}(p_\infty) = 0$, $\nabla W_{g_\infty}(p_\infty) = 0$, and $\nabla^2 W_{g_\infty}(p_\infty) \neq 0$.
\end{itemize}
\end{corollary} 

The question arises as to whether Theorem \ref{blow.up} holds in dimension less than $52$. Recent work of the author and F.~Marques \cite{Brendle-Marques} indicates that similar blow-up examples exist for all $n \geq 25$. Hence, the condition $n \leq 24$ in Theorem \ref{compactness.24} is sharp.

D.~Pollack \cite{Pollack} has shown that every Riemannian metric with positive Yamabe constant is $C^0$ close to a smooth metric for which the Yamabe PDE has an arbitrarily large number of solutions. The solutions constructed in \cite{Pollack} have high energy and high Morse index.

Using the techniques developed in \cite{Brendle4}, one can show that, for $n \geq 52$, every Riemannian metric with positive Yamabe constant can be approximated in $W^{2,s}$ by a metric which allows blow-up: 

\begin{proposition}
\label{density}
Let $M$ be a compact manifold of dimension $n \geq 52$, and let $g_0$ be a smooth Riemannian metric on $M$ with $Y(M,g_0) > 0$. Moreover, let $\beta$ and $s$ be real numbers such that $\beta > 0$ and $s > \frac{n}{2}$. Then we can find a smooth Riemannian metric $g$ on $M$ with the following properties: 
\begin{itemize}
\item[(i)] $\|g - g_0\|_{W^{2,s}(M,g_0)} < \beta$ 
\item[(ii)] There exists a sequence of positive functions $v_\nu \in C^\infty(S^n)$ ($\nu \in \mathbb{N}$) such that 
$\frac{4(n-1)}{n-2} \, \Delta_g v_\nu - R_g \, v_\nu + 4n(n-1) \, v_\nu^{\frac{n+2}{n-2}} = 0$ for all $\nu \in \mathbb{N}$, $E_g(v_\nu) < Y(S^n)$ for all $\nu \in \mathbb{N}$, and $\sup_M v_\nu \to \infty$ as $\nu \to \infty$.
\end{itemize}
\end{proposition}

By the results of Li and Zhang \cite{Li-Zhang1}, this statement is no longer true if we replace the $W^{2,s}$-norm by the $C^2$-norm. 

In order to prove Proposition \ref{density}, we pick an arbitrary point $p \in M$. In geodesic normal 
coordinates around $p$, the metric $g_0$ satisfies $g_{0,ik}(x) = \delta_{ik} + O(|x|^2)$. Let $\eta: \mathbb{R} \to \mathbb{R}$ 
be a smooth cut-off function such that $\eta(t) = 1$ for $t \leq 1$ and $\eta(t) = 0$ for $t \geq 2$. We 
define a trace-free symmetric two-tensor $h(x)$ by 
\[h_{ik}(x) = \sum_{N=N_0}^\infty \eta(4N^2 \, |x - y_N|) \, 2^{-N} \, (2^{-N} - |x - y_N|^2) \, H_{ik}(x - y_N),\] 
where $y_N = (\frac{1}{N},0,\hdots,0) \in \mathbb{R}^n$. Clearly, $h(x)$ is smooth. 
We next define a Riemannian metric $g$ by 
\[g(x) = [1 - \eta(N_0 \, |x|)] \, g_0(x) + \eta(N_0 \, |x|) \, \exp(h(x)).\] 
It is easy to see that 
\[\|g - g_0\|_{W^{2,s}(M,g_0)} < \beta\] 
if $N_0$ is sufficiently large. Using the arguments in \cite{Brendle4}, one can 
construct a sequence of functions $v_\nu$ with the desired properties.

\section{Deformation of Riemannian metrics by their scalar curvature}

We next discuss a natural parabolic PDE associated with the Yamabe problem. Let $M$ be a compact 
manifold of dimension $n \geq 3$, and let $g_0$ be a Riemannian metric on $M$. We consider the 
following flow of Riemannian metrics: 
\begin{equation} 
\label{yamabe.flow.1} 
\frac{\partial}{\partial t} g(t) = -(R_{g(t)} - r_{g(t)}) \, g(t), \qquad g(0) = g_0. 
\end{equation}
As usual, $R_{g(t)}$ denotes the scalar curvature of $g(t)$. Moreover, 
\begin{equation} 
r_{g(t)} = \frac{\int_M R_{g(t)} \, dvol_{g(t)}}{\text{\rm Vol}_{g(t)}(M)} 
\end{equation} 
denotes the mean value of the scalar curvature of $g(t)$. (This guarantees that the volume of $M$ is constant 
along the flow.) The evolution equation (\ref{yamabe.flow.1}) is equivalent to a parabolic PDE for the conformal factor. To see this, we write $g(t) 
= u(t)^{\frac{4}{n-2}} \, g_0$, where $u(t)$ is a positive function defined on $M$. 
If $g(t)$ is a solution of (\ref{yamabe.flow.1}), then $u(t)$ satisfies the PDE 
\begin{equation} 
\label{yamabe.flow.2} 
\frac{\partial}{\partial t} u(t)^{\frac{n+2}{n-2}} = \frac{n+2}{4} \, \Big ( \frac{4(n-1)}{n-2} \, 
\Delta_{g_0} u(t) - R_{g_0} \, u(t) + r_{g(t)} \, u(t)^{\frac{n+2}{n-2}} \Big ). 
\end{equation}
This equation is the gradient flow to the Yamabe functional $E_{g_0}(u)$. It can be 
viewed as the parabolic analogue of the Yamabe PDE (\ref{yamabe.pde}). For that reason, 
(\ref{yamabe.flow.1}) is called the Yamabe flow. 

The Yamabe flow was first proposed by R.~Hamilton. Hamilton showed that the Yamabe flow exists 
for all time (i.e. the flow never develops a singularity in finite time). 

Once longtime existence is known, the next issue is the asymptotic behavior of the flow as $t \to \infty$. 
This is not a difficult problem if $Y(M,g_0) \leq 0$: in this case, the flow converges exponentially to a 
metric of constant scalar curvature. However, the problem is much more subtle if $Y(M,g_0) > 0$. 
B.~Chow \cite{Chow} proved that the Yamabe flow converges to a metric of constant scalar 
curvature provided that $(M,g_0)$ is locally conformally flat and $g_0$ has positive Ricci 
curvature. The second assumption was subsequently removed by R.~Ye \cite{Ye}: 
Ye showed that the flow converges to a metric of constant scalar curvature whenever $(M,g_0)$ 
is locally conformally flat. The proof in \cite{Ye} rests on the injectivity of the developing map 
and a reflection argument. This method uses in a fundamental way the assumption that $(M,g_0)$ is 
locally conformally flat. 

H.~Schwetlick and M.~Struwe \cite{Schwetlick-Struwe} introduced a new approach based on 
concentration-compactness arguments. In low dimensions, they were able to prove 
the convergence of the Yamabe flow assuming only a bound on the initial energy. More precisely, 
if $3 \leq n \leq 5$ and the Yamabe energy of the initial metric is less than 
$\big [ Y(M,g_0)^{\frac{n}{2}} + Y(S^n)^{\frac{n}{2}} \big ]^{\frac{2}{n}}$, then the Yamabe flow 
converges to a metric of constant scalar curvature as $t \to \infty$. In \cite{Brendle1}, 
we proved the convergence of the Yamabe flow for initial metrics with arbitrary energy: 

\begin{theorem}
\label{convergence.of.flow.1}
Let $(M,g_0)$ be a compact Riemannian manifold of dimension $n$. If $3 \leq n \leq 5$, then the 
Yamabe flow converges to a metric of constant scalar curvature as $t \to \infty$.
\end{theorem}

In the remainder of this section, we will sketch the main steps involved in the proof of Theorem \ref{convergence.of.flow.1}. The argument is particularly simple if $(M,g_0)$ is the round sphere; see \cite{Brendle2}.

Since the volume of $M$ is preserved by the flow, we may assume that $\text{\rm Vol}_{g(t)}(M) = 1$ 
for all $t \geq 0$. For abbreviation, let $r_\infty = \lim_{t \to \infty} r_{g(t)}$. (The limit 
exists since the function $t \mapsto r_{g(t)}$ is decreasing.) It follows from work of Schwetlick and Struwe \cite{Schwetlick-Struwe} that 
\[\lim_{t \to \infty} \int_M |R_{g(t)} - r_\infty|^s \, dvol_{g(t)} = \lim_{t \to \infty} 
\int_M |R_{g(t)} - r_{g(t)}|^s \, dvol_{g(t)} = 0\] for $1 < s < \frac{n+2}{2}$. Hence, if 
$t_\nu$ ($\nu \in \mathbb{N}$) is a sequence of times such that $t_\nu \to \infty$ as 
$\nu \to \infty$, then the functions $u(t_\nu)$ form a Palais-Smale sequence. The 
properties of such sequences are described by the following important result:

\begin{proposition} 
\label{blow.up.analysis}
Let $c$ be a positive real number, and let $u_\nu$ ($\nu \in \mathbb{N}$) be a sequence of positive functions 
satisfying 
\[\int_M u_\nu^{\frac{2n}{n-2}} \, dvol_{g_0} = 1\] 
for all $\nu \in \mathbb{N}$ and 
\[\frac{4(n-1)}{n-2} \, \Delta_{g_0} u_\nu - R_{g_0} \, u_\nu + c \, u_\nu^{\frac{n+2}{n-2}} 
\to 0\] 
in $L^{\frac{2n}{n+2}}(M)$. After passing to a subsequence if necessary, we can find a 
non-negative integer $m$, a family of points $p_{k,\nu} \in M$ ($1 \leq k \leq m$, $\nu \in \mathbb{N}$), and 
a family of positive real numbers $\varepsilon_{k,\nu}$ ($1 \leq k \leq m$, $\nu \in \mathbb{N}$) 
with the following properties: 
\begin{itemize}
\item[(i)] For each $k$, we have $\varepsilon_{k,\nu} \to 0$ as $\nu \to \infty$. 
\item[(ii)] For all $i \neq j$, we have 
\[\frac{\varepsilon_{i,\nu}}{\varepsilon_{j,\nu}} + \frac{\varepsilon_{j,\nu}}{\varepsilon_{i,\nu}} 
+ \frac{d(p_{i,\nu},p_{j,\nu})^2}{\varepsilon_{i,\nu} \, \varepsilon_{j,\nu}} \to \infty\] 
as $\nu \to \infty$.
\item[(iii)] There exists a non-negative function $u_\infty$ such that 
\[u_\nu - \sum_{k=1}^m \Big ( \frac{4n(n-1)}{c} \Big )^{\frac{n-2}{4}} \, \Big ( \frac{\varepsilon_{k,\nu}}{\varepsilon_{k,\nu}^2 + d(p_{k,\nu},\cdot)^2} \Big )^{\frac{n-2}{2}} 
\to u_\infty\] 
in $W^{1,2}(M,g_0)$. 
\item[(iv)] The function $u_\infty$ is a smooth solution of the PDE 
\[\frac{4(n-1)}{n-2} \, \Delta_{g_0} u_\infty - R_{g_0} \, u_\infty + c \, u_\infty^{\frac{n+2}{n-2}} = 0.\] 
\item[(v)] If $u_\infty$ does not vanish identically, then $u_\infty$ is positive everywhere, 
and we have $c = \big [ E_{g_0}(u_\infty)^{\frac{n}{2}} + m \, Y(S^n)^{\frac{n}{2}} \big ]^{\frac{2}{n}}$. 
Otherwise, we have $c = \big [ m \, Y(S^n)^{\frac{n}{2}} \big ]^{\frac{2}{n}}$. 
\end{itemize}
\end{proposition}

Proposition \ref{blow.up.analysis} is a variant of a result of M.~Struwe \cite{Struwe1}. 
Property (ii) is due to A.~Bahri and J.M.~Coron \cite{Bahri-Coron}. 

Our goal is to prove that the Yamabe flow converges to a metric of constant
scalar curvature as $t \to \infty$. To that end, we need to show that
volume concentration does not occur (i.e. $m = 0$). The following result
is the key ingredient in the proof of Theorem \ref{convergence.of.flow.1}:

\begin{proposition}
\label{lojasiewicz.inequality}
Assume that $3 \leq n \leq 5$. Let $c$ be a positive real number, and let $u_\nu$ ($\nu \in \mathbb{N}$) be a sequence of positive functions 
satisfying 
\[\int_M u_\nu^{\frac{2n}{n-2}} \, dvol_{g_0} = 1\] 
for all $\nu \in \mathbb{N}$ and 
\[\frac{4(n-1)}{n-2} \, \Delta_{g_0} u_\nu - R_{g_0} \, u_\nu + c \, u_\nu^{\frac{n+2}{n-2}} 
\to 0\] 
in $L^{\frac{2n}{n+2}}(M)$. Then we can find positive real numbers $\gamma$ and $C$ such that 
\begin{align*} 
&E_{g_0}(u_\nu) - c \\ 
&\leq C \, \bigg ( \int_M \Big | \frac{4(n-1)}{n-2} \, \Delta_{g_0} u_\nu - R_{g_0} \, u_\nu + c \, u_\nu^{\frac{n+2}{n-2}} 
\Big |^{\frac{2n}{n+2}} \, dvol_{g_0} \bigg )^{\frac{n+2}{2n} \, (1+\gamma)} 
\end{align*}
for infinitely many $\nu \in \mathbb{N}$.
\end{proposition}

If $u_\infty$ vanishes, then we can choose $\gamma = 1$. The same is true if $u_\infty > 0$ 
is a non-degenerate solution of the Yamabe PDE. 

In \cite{Brendle1}, Proposition \ref{lojasiewicz.inequality} was stated for 
sequences of the form $u_\nu = u(t_\nu)$, where $g(t) = u(t)^{\frac{4}{n-2}} \, g_0$ 
is a solution of the Yamabe flow (see \cite{Brendle1}, Proposition 3.3). However, the proof works 
for general Palais-Smale sequences. The following result is an immediate consequence of 
Proposition \ref{lojasiewicz.inequality}: 

\begin{corollary} 
\label{lojasiewicz.inequality.flow.version}
Assume that $3 \leq n \leq 5$. Suppose that $g(t) = u(t)^{\frac{4}{n-2}} \, g_0$ is a solution to the Yamabe flow 
such that 
\[\int_M u(t)^{\frac{2n}{n-2}} \, dvol_{g_0} = 1\] 
for all $t \geq 0$. Let $r_\infty = \lim_{t \to \infty} r_{g(t)}$. 
Then there exist positive real numbers $\gamma$ and $t_0$ such that 
\begin{align*} 
&E_{g_0}(u(t)) - r_\infty \\ 
&\leq \bigg ( \int_M \Big | \frac{4(n-1)}{n-2} \, 
\Delta_{g_0} u(t) - R_{g_0} \, u(t) + r_\infty \, u(t)^{\frac{n+2}{n-2}} 
\Big |^{\frac{2n}{n+2}} \, dvol_{g_0} \bigg )^{\frac{n+2}{2n} \, (1+\gamma)} 
\end{align*}
for all $t \geq t_0$.
\end{corollary}

Corollary \ref{lojasiewicz.inequality} can be used to rule out volume concentration along 
the Yamabe flow. Indeed, it follows from Corollary \ref{lojasiewicz.inequality.flow.version} 
and (\ref{transformation.law.for.scalar.curvature}) that 
\[r_{g(t)} - r_\infty \leq \bigg ( \int_M |R_{g(t)} - r_\infty|^{\frac{2n}{n+2}} \, dvol_{g(t)} \bigg )^{\frac{n+2}{2n} \, (1+\gamma)}\] 
for $t \geq t_0$. This implies 
\[r_{g(t)} - r_\infty \leq C \, \bigg ( \int_M |R_{g(t)} - r_{g(t)}|^{\frac{2n}{n+2}} \, dvol_{g(t)} \bigg )^{\frac{n+2}{2n} \, (1+\gamma)}\] 
if $t$ is sufficiently large. On the other hand, we have 
\[r_{g(t)} - r_\infty = \frac{n-2}{2} \int_t^\infty \int_M (R_{g(\tau)} - r_{g(\tau)})^2 \, dvol_{g(\tau)} \, d\tau\] 
for all $t \geq 0$. Putting these facts together, one can show that 
\[\int_0^\infty \bigg ( \int_M (R_{g(\tau)} - r_{g(\tau)})^2 \, dvol_{g(\tau)} \bigg )^{\frac{1}{2}} \, d\tau < \infty.\] 
Hence, given any $\eta > 0$, we can find a positive real number $r$ such that 
\[\int_{B_r(p)} u(t)^{\frac{2n}{n-2}} \, dvol_{g_0} \leq \eta\] 
for all $p \in M$ and all $t \geq 0$. (Here, $B_r(p)$ is a geodesic ball 
with respect to the background metric $g_0$.) This shows that volume concentration 
does not occur. Hence, the flow converges to a metric of constant scalar curvature 
as $t \to \infty$.

Therefore, it suffices to prove Proposition \ref{lojasiewicz.inequality}. 
The proof of this result uses a well-known inequality for analytic functions 
due to Lojasiewicz (cf. \cite{Simon}). More importantly, we need a family of functions 
$\overline{u}_{(p,\varepsilon)}$ ($p \in M$, $0 < \varepsilon \leq \varepsilon_0$) 
with the following properties:
\begin{itemize}
\item The Yamabe energy of $\overline{u}_{(p,\varepsilon)}$ is bounded by 
\[E_{g_0}(\overline{u}_{(p,\varepsilon)}) \leq Y(S^n)\] 
for all $p \in M$ and $0 < \varepsilon \leq \varepsilon_0$. 
\item As $\varepsilon \to 0$, the functions $\overline{u}_{(p,\varepsilon)}$ 
converge to the "standard bubble" on $\mathbb{R}^n$ after rescaling. More precisely, 
\[\lim_{\varepsilon \to 0} \varepsilon^{\frac{n-2}{2}} \, 
\overline{u}_{(\varepsilon,p)}(\exp_p(\varepsilon \xi)) = \Big ( \frac{1}{1 + |\xi|^2} \Big )^{\frac{n-2}{2}}\] 
for all $p \in M$ and $\xi \in T_p M$.
\end{itemize}
For $3 \leq n \leq 5$, the existence of a family of functions with these properties 
follows work of R.~Schoen \cite{Schoen1}.

\section{Convergence of the Yamabe flow for $n \geq 6$}

In this section, we provide sufficient conditions for the convergence of the Yamabe flow 
in higher dimensions. Throughout this section, we fix a compact Riemannian manifold $(M,g)$ of dimension 
$n \geq 6$ with positive Yamabe constant.

\begin{definition} 
Let $d = [\frac{n-2}{2}]$. We denote by $\mathcal{Z}$ the set of all points $p \in M$ such that 
\[\limsup_{q \to p} d(p,q)^{2-d} \, |W_g(q)| = 0,\] 
where $W_g$ denotes the Weyl tensor of $g$.
\end{definition}

In other words, a point $p \in M$ belongs to $\mathcal{Z}$ if and only if the Weyl tensor 
vanishes to an order greater than $[\frac{n-6}{2}]$ at $p$. Clearly, $\mathcal{Z}$ is a 
closed subset of $M$. Moreover, the set $\mathcal{Z}$ depends only on the conformal class 
of $g$.

\begin{theorem}
\label{convergence.of.flow.2}
Assume that $\mathcal{Z} = \emptyset$. Then, for every initial metric in the conformal class of $g$, 
the Yamabe flow converges to a metric of constant scalar curvature as $t \to \infty$.
\end{theorem}

The condition $\mathcal{Z} = \emptyset$ can be removed if the
positive mass theorem holds. To explain this, consider a point
$p \in \mathcal{Z}$, and let $G$ be the Green's function of the 
conformal Laplacian with pole at $p$. Then the manifold $(M \setminus \{p\},
G^{\frac{4}{n-2}} \, g)$ is scalar flat and satisfies the decay 
conditions of R.~Bartnik \cite{Bartnik}. Hence, this manifold has a well-defined ADM mass. 

We can prove the convergence of the Yamabe flow provided that the ADM 
mass of $(M \setminus \{p\},G^{\frac{4}{n-2}} \, g)$ is positive for 
all $p \in \mathcal{Z}$. If $M \setminus \{p\}$ is a spin manifold, 
then the positivity of the ADM mass follows from work of Bartnik 
(see \cite{Bartnik}, Theorem 6.3; see also \cite{Schoen-Yau1},\cite{Witten}).

\begin{theorem}
\label{convergence.of.flow.3}
Assume that $M \setminus \{p\}$ is a spin manifold. Then, for every choice of 
the initial metric $g_0$, the Yamabe flow converges
to a metric of constant scalar curvature as $t \to \infty$.
\end{theorem}

In order to prove Theorem \ref{convergence.of.flow.2}, we need to construct 
a suitable family of test functions $\overline{u}_{(p,\varepsilon)}$. If the 
Weyl tensor of $g$ is non-zero at each point on $M$, we can use a result 
of T.~Aubin \cite{Aubin1} to construct a family of test functions with the required 
properties. In the remaining part of this section, we sketch how this construction can be generalized to deal with points where the Weyl tensor is equal to $0$. The details can be found in \cite{Brendle3}.

Consider a point $p \in M$. After changing the metric $g$ conformally, we may assume that 
$\det g(x) = 1 + O(|x|^{2d+2})$ in geodesic normal coordinates around $p$ 
(see \cite{Lee-Parker}, Theorem 5.1, or \cite{Schoen-Yau2}, Theorem 3.1). Hence, we may write 
$g(x) = \exp(h(x))$, where $h(x)$ is a symmetric two-tensor satisfying 
$\text{\rm tr} \, h(x) = O(|x|^{2d+2})$ and $\sum_{k=1}^n h_{ik}(x) \, x_k = 0$ 
for $i = 1, \hdots, n$. Let $H_{ik}(x)$ be the Taylor polynomial of order $d$ associated with 
the function $h_{ik}(x)$. Clearly, $H_{ik}(x)$ is a trace-free symmetric two-tensor on 
$\mathbb{R}^n$ satisfying $\sum_{k=1}^n H_{ik}(x) \, x_k = 0$ for $i = 1, \hdots, n$. 
Note that $H_{ik}(x) \not\equiv 0$ since $\mathcal{Z} = \emptyset$.

We denote by $G$ the Green's function of the conformal Laplacian with pole at $p$. Let 
$\chi: \mathbb{R} \to \mathbb{R}$ is a non-negative cutoff function satisfying $\chi(t) = 1$ for $t \leq \frac{4}{3}$ and $\chi(t) = 0$ for $t \geq \frac{5}{3}$. 
Moreover, suppose that $\delta$ and $\varepsilon$ are positive real numbers such that $\delta \geq 2\varepsilon$. 
We define a function $u_\varepsilon$ by 
\[u_\varepsilon(x) = \Big ( \frac{\varepsilon}{\varepsilon^2 + |x|^2} \Big )^{\frac{n-2}{2}}.\] 
We can find a vector field $V^{(\varepsilon,\delta)}$ on $\mathbb{R}^n$ such that 
\[\sum_{k=1}^n \partial_k \Big [ u_\varepsilon^{\frac{2n}{n-2}} \, (\chi(|x|/\delta) \, H_{ik} - \partial_i V_k^{(\varepsilon,\delta)} - \partial_k V_i^{(\varepsilon,\delta)} + \frac{2}{n} \, \text{\rm div} \, V^{(\varepsilon,\delta)} \, \delta_{ik}) \Big ] = 0\] 
for $i = 1, \hdots, n$. We then define a function $v_{(\varepsilon,\delta)}: M \to \mathbb{R}$ by 
\[v_{(\varepsilon,\delta)} = \chi(|x|/\delta) \, (u_\varepsilon + w_{(\varepsilon,\delta)}) + (1 - \chi(|x|/\delta)) \, \varepsilon^{\frac{n-2}{2}} \, G,\] 
where 
\[w_{(\varepsilon,\delta)} = \sum_{l=1}^n \partial_l u_\varepsilon \, V_l^{(\varepsilon,\delta)} + \frac{n-2}{2n} \, u_\varepsilon
\, \text{\rm div} \, V^{(\varepsilon,\delta)}.\] 
Note that $w_{(\varepsilon,\delta)}$ is a solution of the linear elliptic PDE 
\[\Delta w_{(\varepsilon,\delta)} + n(n+2) \, u_\varepsilon^{\frac{4}{n-2}} \, w_{(\varepsilon,\delta)} = \frac{n-2}{4(n-1)} \, u_\varepsilon \, \sum_{i,k=1}^n \partial_i \partial_k \big [ \chi(|x|/\delta) \, H_{ik} \big ].\] 
The function $w_{(\varepsilon,\delta)}$ can be viewed as a correction term which is needed to 
push the Yamabe energy of $v_{(\varepsilon,\delta)}$ below $Y(S^n)$. The idea of adding a small 
correction term to $u_\varepsilon$ goes back to work of E.~Hebey and M.~Vaugon \cite{Hebey-Vaugon}. 

Our goal is to show that $E_g(v_{(\varepsilon,\delta)}) < Y(S^n)$ if $\varepsilon$ is sufficiently small. 
The leading term in the asymptotic expansion for 
\[\Big ( \frac{Y(S^n)}{4n(n-1)} \Big )^{\frac{n-2}{2}} \, \big ( E_g(v_{(\varepsilon,\delta)}) - Y(S^n) \big )\] 
is given by 
\begin{align*} 
&- \int_{B_\delta(0)} \frac{1}{4} \, u_\varepsilon^2 \, \sum_{i,k,l=1}^n \partial_l H_{ik} \, \partial_l H_{ik} + \int_{B_\delta(0)} \frac{1}{2} \, u_\varepsilon^2 \, \sum_{i,k,l=1}^n \partial_k H_{ik} \, \partial_l H_{il} \\ 
&+ \int_{B_\delta(0)} 2 \, u_\varepsilon \, w_{(\varepsilon,\delta)} \, \sum_{i,k=1}^n \partial_i \partial_k H_{ik} + \int_{B_\delta(0)} \frac{4(n-1)}{n-2} \, |dw_{(\varepsilon,\delta)}|^2 \\ 
&- \int_{B_\delta(0)} \frac{4(n-1)}{n-2} \, n(n+2) \, u_\varepsilon^{\frac{4}{n-2}} \, w_{(\varepsilon,\delta)}^2. 
\end{align*} 
In order to estimate this term, we use the identity 
\begin{align} 
\label{key.identity}
&- \frac{1}{4} \, u_\varepsilon^2 \, \sum_{i,k,l=1}^n \partial_l H_{ik} \, \partial_l H_{ik} + \frac{1}{2} \, u_\varepsilon^2 \, \sum_{i,k,l=1}^n \partial_k H_{ik} \, \partial_l H_{il} \notag \\ 
&+ 2 \, u_\varepsilon \, w_{(\varepsilon,\delta)} \, \sum_{i,k=1}^n \partial_i \partial_k H_{ik} + \frac{4(n-1)}{n-2} \, |dw_{(\varepsilon,\delta)}|^2 \notag \\ 
&- \frac{4(n-1)}{n-2} \, n(n+2) \, u_\varepsilon^{\frac{4}{n-2}} \, w_{(\varepsilon,\delta)}^2 \notag \\ 
&= - \frac{1}{4} \sum_{i,k,l=1}^n |Q_{ik,l}^{(\varepsilon,\delta)}|^2 - 2 \, u_\varepsilon^{\frac{2n}{n-2}} \, \sum_{i,k=1}^n |T_{ik}^{(\varepsilon,\delta)}|^2 \\ 
&+ \frac{1}{2} \, \sum_{i=1}^n \bigg | \sum_{k=1}^n \Big ( u_\varepsilon \, \partial_k T_{ik}^{(\varepsilon,\delta)} + \frac{2n}{n-2} \, \partial_k u_\varepsilon \, T_{ik}^{(\varepsilon,\delta)} \Big ) \bigg |^2 \notag \\ 
&+ \text{\rm divergence terms}. \notag
\end{align} 
Here, $T_{ik}^{(\varepsilon,\delta)}$ and $Q_{ik,l}^{(\varepsilon,\delta)}$ are defined by 
\[T_{ik}^{(\varepsilon,\delta)} = H_{ik} - \partial_i V_k^{(\varepsilon,\delta)} - \partial_k V_i^{(\varepsilon,\delta)} + \frac{2}{n} \, \text{\rm div} \, V^{(\varepsilon,\delta)} \, \delta_{ik}\] 
and 
\begin{align*}
Q_{ik,l}^{(\varepsilon,\delta)} 
&= u_\varepsilon \, \partial_l T_{ik}^{(\varepsilon,\delta)} - \frac{2}{n-2} \, \partial_i u_\varepsilon \, T_{kl}^{(\varepsilon,\delta)} - \frac{2}{n-2} \, \partial_k u_\varepsilon \, T_{il}^{(\varepsilon,\delta)} \\
&+ \frac{2}{n-2} \, \partial_p u_\varepsilon \, T_{ip}^{(\varepsilon,\delta)} \, \delta_{kl} + \frac{2}{n-2} \, \partial_p u_\varepsilon \, T_{kp}^{(\varepsilon,\delta)} \, \delta_{il}. 
\end{align*}
The identity (\ref{key.identity}) is closely related to the second variation formula for the Einstein-Hilbert action near the round metric on $S^n$ (see \cite{Besse}, Section 4G). 
The first two terms on the right hand side are non-positive. Moreover, the third term vanishes 
in $B_\delta(0)$. To see this, observe that 
\[\sum_{k=1}^n \partial_k \big ( u_\varepsilon^{\frac{2n}{n-2}} \, T_{ik}^{(\varepsilon,\delta)} \big ) = 0\] 
for $|x| \leq \delta$. This implies 
\[\sum_{k=1}^n \Big ( u_\varepsilon \, \partial_k T_{ik}^{(\varepsilon,\delta)} + \frac{2n}{n-2} \, 
\partial_k u_\varepsilon \, T_{ik}^{(\varepsilon,\delta)} \Big ) = 0\] 
for $|x| \leq \delta$.

\end{document}